\documentclass[12pt]{article}
\usepackage{graphics}
\usepackage{color}
\usepackage{epsfig}
\usepackage{latexsym}
\usepackage{amssymb}
\usepackage{multicol}
\usepackage{graphicx}
\usepackage{amsmath}
\usepackage{amsfonts}
\DeclareMathSymbol\square {\mathrel}{AMSa}{"03}
\newtheorem{theorem}{Theorem }
\newtheorem{cor}{Corollary}

\title {Heron's Formula, Descartes Circles, and Pythagorean Triangles}
\author{Frank Bernhart and H. Lee Price}

\begin{document} 
\date{January 1, 2007}

\maketitle



\begin{abstract}
This article highlights interactions of diverse areas: the Heron 
formula for the area of a triangle, the Descartes circle equation, and right 
triangles with integer or rational sides. New and old results are 
synthesized. We first exploit elementary observations about circles to characterize an 
arbitrary triangle using three circles. A fourth circle brings a certain 
symmetry -- the four radii are exactly the factors of the Heron area 
formula. The four equi-circles  also play a role. 

When the triangle is a right triangle, the four circles combine into a 
single tangent cluster. Their centers form a rectangle, the system is 
reflection-congruent to its own dual set, and the dual circles correspond to 
the equi-circles (in-circle and ex-circles). The four radii satisfy 
$r_4=r_1 +r_2 +r_3, \ r_2 \cdot r_3= r_1 \cdot r_4$. Moreover, $[-r_1 ,r_2 ,r_3 ,r_3 ]$ 
is a quadruple that satisfies the Descartes Circle Equation \textbf{(DCE)}: the sum of the 
squares is half the square of the sum. Thus a right triangle with legs 
$2=r_1 + r_2, \ 4=r_1 +r_3$ and hypotenuse $2\sqrt{5}=r_2 +r_3$ gives this \textbf{DCE} solution:
$$ [\sqrt 5 -  3, \ \sqrt 5  -  1, \  \sqrt 5  +  1, \  \sqrt 5  +  3]$$
Pythagorean Triples (\textbf{PT}'s) are positive integers that form the legs 
$a,b$ and the hypotenuse ${c}$ of a right triangle. The basic theory fits quite well with 
the circle method. Convenient radial equation $2{r_1}^{ 2} =  (r_2 -  r_1) (r_3 -  r_1)$ 
is equivalent to an old and powerful identity of Dickson:
$$ \tfrac{1}{2}(a+b+c)^2= (c-b)(c-a)$$

Ratios $\tfrac{r_1}{r_2}, \ \tfrac{r_1}{r_3}$ are the half-angle tangents for the right triangle; the reduced numerators and denominators form a Fibonacci-rule sequence $\cal P$. Primitive sequences 
$\cal P$ and their multiples serve to generate primitive \textbf{PT}'s and their multiples. This includes 
the ancient standard solutions, and new ones. The product of the four numbers of $\cal P$ is the area $\textit{G}$ of the triangle.

One can replace the in-radius $r_1$ with any ex-radius ($r_2,$ $r_3,$ or $r_4$) by 
rearranging and adjusting $\cal P$. This provides a new and natural construction of the 
Barning-Hall ternary tree, which grows all primitive \textbf{PT}'s from the root triple $[3,4,5]$. 
A classical theorem involving the nine-point circle provides yet another tree construction. Every primitive \textbf{PT} when rescaled by factor $1/\textit{G}$ furnishes a Descartes quadruple of tangent 
circles with integral curvatures, and so generates an integral Apollonian packing \textbf{(IAP)} containing a 
rectangle of centers. Thus \textbf{PT}'s serve to generate an infinite number of in-equivalent integral packings. 
A table of quadruples ends the article.
\end{abstract}

\section{Going in Circles.}

A plane circle $ {\cal C}( {\rm A}, r)\ $ has center ${\rm A}$ and radius $r > 0$. In \textbf{Fig.1(a)} 
we find three circles of arbitrary radii $r_1 ,r_2 ,r_3 $ arranged with 
their centers on a line. Each one is (externally) tangent to its neighbors. A final circle is added, its center on the same line, tangent as shown 
(internally), and with radius $r_4 = r_1 + r_2 + r_3$. 

\begin{figure}[htbp]
\centerline{\includegraphics[width=2.75in]{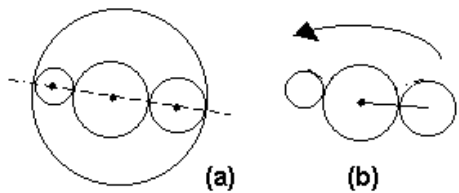}}
\caption{Three circles and their `\! sum \!'.}
\end{figure}

A triad of three tangent circles is now constructed, in four different ways. 
First omit one circle, leaving two circles tangent to the third. One of the 
two is then rolled around the third (externally or internally) as in 
\textbf{Fig.1(b)} until all three are tangent. In \textbf{Fig.2} there is one such triad. 

\begin{figure}[htbp]
\centerline{\includegraphics[width=2.6in]{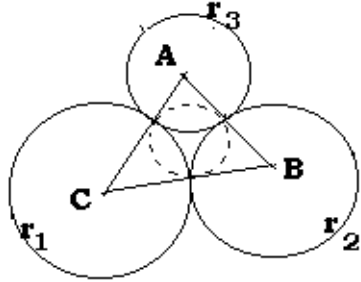}}
\caption {A tangent triad defines a triangle.}
\end{figure}

The centers $ \rm A,B,C$ are connected to make a triangle, with side 
$ a$ opposite vertex $ \rm A$, etc. We easily find (including $r_4$ for later convenience):

\begin{equation}
	\begin{array}{lll}
  a  = & r_1 \ + r_2 \ =  &  r_4 - \ r_3  \\ 
  b  = & r_1 \ + r_3 \ =  &  r_4 - \ r_2  \\ 
  c  = & r_2 \ + r_3 \ =  &  r_4 - \ r_1  
 \end{array}
\end{equation}

Can we guarantee that the sum of two sides always is greater than the third? 
That $ \ a+b-c \ $ for example is positive? Investigating, we are pleased to find  \textbf{(2)}.
\begin{equation}
	\begin{array}{rl}
	 r_1 = & \tfrac{1}{2} (+ a + b - c)  \\
	 r_2 = & \tfrac{1}{2} (+ a - b + c)  \\ 
	 r_3 = & \tfrac{1}{2} (- a + b + c)  \\ 
	 r_4 = & \tfrac{1}{2} (+ a + b + c)  \\
	 {} =  & (r_1 + r_2 + r_3)
 \end{array}
\end{equation}

In other words, the three diameters $2r_1  ,\ 2r_2  , \ 2r_3$ are precisely those 
values which must be positive in order for $ a, b, c$ to form a non-degenerate triangle. 
As a bonus, the value $ r_4 \ $ is seen to be half of the perimeter!

Before inspecting the remaining triads, we point out in \textbf{Fig.2} 
the in-circle of the triangle. It is well known and easy to show that the 
in-circle is orthogonal to all three circles. It is also the circle 
determined by the contact points of the triad.

Now for the other triads. Oddly enough, joining the centers gives exactly 
the same triangle! \textbf{Fig.3} shows that the ex-circles play the same 
role, \textit{mutatis mutandis}, that the in-circle played in the previous figure.

\begin{figure}[htbp]
\centerline{\includegraphics[width=3.0in]{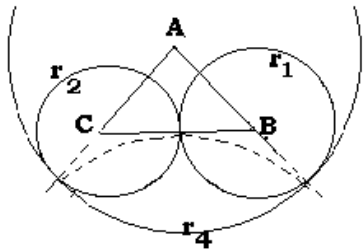}}
\caption {Another way to obtain the same triangle.}
\end{figure}

Until now, the fourth circle seemed somewhat unnecessary, but now it is seen 
as essential if the four triads are to form an integrated group. The central 
fact that all four triads give the same triangle is expressed by the easy 
equations \textbf{(1), (2)}.

Let ${\cal R}{ =  [ r_1,     r_ 2,  r_3,   r_4   ]}$ be positive parameters, arbitrary 
except for $r_4 ,$ as given. Then $\cal R$ makes a rather interesting description of the shape 
of an arbitrary triangle. The triangle area $\textit{G}$ is neatly expressed as the positive 
square root of the product of all four parameters.

\begin{equation}
{\textit{G}}^{2} =  \prod   r_i  =  r_1 \cdot  r_2  \cdot r_3 \cdot r_4	
\end{equation}

This is a ``circular'' equivalent to a classical formula, known as Heron's 
formula. In the usual way of expressing this result, we have semi-perimeter $ s $ and factors $s - a, \ s - b, \ s - c{\  }$. However, the ``translation'' provided by \textbf{(4)} is easy to check.
\begin{equation}
	\begin{array}{rl}
		s    = & r_4  = {\tfrac{1}{2}} (+a+b+c) \\ 
		s-a  = & r_3  = {\tfrac{1}{2}} (-a+b+c) \\ 
		s-b  = & r_2  = {\tfrac{1}{2}} (+a-b+c) \\ 
		s-c  = & r_1  = {\tfrac{1}{2}} (+a+b-c) \\ 
 \end{array}
\end{equation}

Our presentation has the advantage of using ``circular reasoning'' to give 
each factor on the right-hand-side of {(3)} a basic geometric 
meaning. We denote by ${\cal T}_i $ the tangent triad formed by circles with 
radii $r_j ,  j \ne {\ i}$. The three contact points in the triad determine a 
``companion'' orthogonal circle, which is an in-circle ${\cal T}_4$ or an ex-circle 
$({\cal T}_1, {\cal T}_2, {\cal T}_3)$ for the triangle of centers. Let the circle ${\cal T}_i $ have radius 
denoted by $s_i$. The in-circle is ${\cal C}(\rm I, \ s_4 )$, and it has an in-center $\rm I$.

If the center $\text I $ is connected to triangle corners $\rm A, B, C $ three 
sub-triangles $\rm ABI$, $\rm BCI$, $\rm CAI$ result. By adding their areas it is 
readily found that ${\textit{G}} =  r_4 \cdot \ s_4$. Very similar reasoning leads to \textbf{(5)}.
 
\begin{equation}
{\textit{G}} =  r_1  \cdot   s_1 =  r_2 \cdot   s_2   =  r_3 \cdot   s_3  =   r_4  \cdot   s_4 	
\end{equation}
Comparing \textbf{(5)} with \textbf{(3)} we discover that ${\textit{G}}^2  =  \prod s_i$. But much more is to come.

\section{Getting it Right.}

Carbon under great pressure may produce a diamond. In our case, it is pressure enough 
to specialize to right triangles. Let the angle at vertex $\rm C$ be a right angle. In 
other words, the Pythagorean equation \textbf{(6)} holds.

\begin{equation}
a^2 + b^2 = c^2 
\end{equation}
Substituting from \textbf{(1)}, we find that \textbf{(6)} is transformed to $r_2  \cdot r_3 =  r_1  \cdot r_4$, 
hence \textbf{(7)}. Also \textbf{(3)} and \textbf{(5)} are simplified to \textbf{(8)} and \textbf{(9)}.

\begin{equation}
r_4 \  =  r_1 +   r_2  +  r_3 \ , \ r_2 \cdot r_3 \  = \ r_1  \cdot r_4 
\end{equation}

\begin{equation}
{\textit{G}} \ = \ r_2  \cdot r_3 \  = \ r_1  \cdot r_4 
\end{equation}

\begin{equation}
r_1 = s_4,  \quad r_2 = s_3,  \quad r_3  = s_2,  \quad r_4  =  s_1 
\end{equation}
 The identity of the $s_i$ and the $r_i$ is 
striking. It correctly forecasts that the case of the right triangle may 
divulge further special features. The first payoff is the coalescing of the 
four triads into a single diagram. For this we need a new section.

\section{Congruent/Dual Tangent Systems.} First we define the four points
$$
\rm C  =  (0, 0), \quad \rm B  =  (a, 0), \quad \rm A  =  (0, b), \quad \rm D  =  (a, b)
$$
These form an $a \times b$ rectangle, as in Fig.4(a), so that 
$$
a = \overline {\rm CB} = \overline {\rm AD}, \quad b  =  \overline {\rm AC} = \overline {\rm BD}, \quad c  =  \overline {\rm AB} = \overline {\rm CD} 
.$$
Here $\rm ABC$ is a right triangle with side $a$ in the X-axis, side $b$ in the Y-axis, and right angle C at the origin.
 
Triangles $\rm ABD$, $\rm ACD$, $\rm BCD$ are congruent copies of the same right triangle. In 
\textbf{Fig.4(a)} below the specific right triangle $[a, b, c]  =  [3, 4, 5]$ is selected for definiteness (and so ${\cal R} = [1, 2, 3, 6  ]$). This diagram can serve for the general case. The four following circles are 
added (Note: $\rm D$ is not usually on circle ${\cal K}_3$). 

\begin{figure}[htbp]
\centerline{\includegraphics[width=3.75in,height=1.66in]{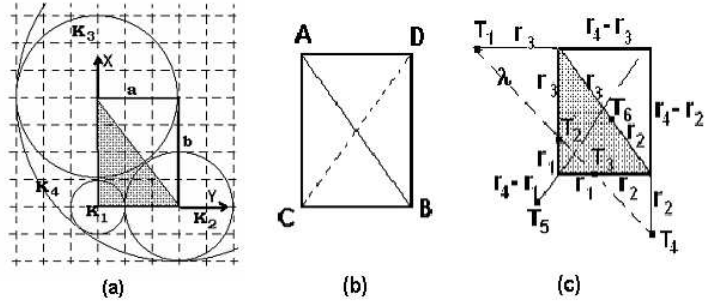}}
\caption{How a rectangle generates four tangent circles.}
\end{figure}
$$
{\cal K}_1 = {\cal C}(\rm C, r_1), \quad  {\cal K}_2 = {\cal C}(\rm B, r_2), \quad {\cal K}_3 = {\cal C }(\rm A, r_3), \quad  {\cal K}_4 = {\cal C}(\rm D, r_4) 
$$
Through \textbf{(1)} we verify easily that circles ${\cal K}_i$
are pairwise tangent. This tangent system will be denoted $\alpha$. 
We see that contact points
$$
{\rm T}_1 =  {\cal K}_3{\cal K}_4, \quad     {\rm T}_2 =  {\cal K}_1{\cal K}_3,   \quad    {\rm T}_3 =  {\cal K}_1{\cal K}_2,   \quad    {\rm T}_4 =  {\cal K}_2{\cal K}_4 
$$
lie on extended sides $\overline {\rm AD}, \  \overline {\rm AC}, \ \overline {\rm CB}, \  \overline {\rm BD}$
, and are given by
$$
{\rm T}_1  = (-r_3, b), \quad {\rm T}_2 =  (0, r_1), \quad {\rm T}_3 = (r_1, 0), \quad  {\rm T}_4 =  (a,  -r_2)
$$

In \textbf{Fig.4(c)} they lie from upper left to lower right on the 
line $\lambda : x  +  y  =  r_1$. The remaining two contacts are 
${\rm T}_5 = {\cal K}_1   {\cal K}_4, \ {\rm T}_6 = {\cal K}_2  {\cal K}_3$. 
The segments $\overline {{\rm C}{{\rm T}_5}}, \ \overline {{\rm A}{{\rm T}_6}}, \ \overline {{\rm B}{{\rm T}_6}},$ are radii of lengths $r_1, \ r_3, \ r_2$. They are also hypotenuses of triangles similar to $\rm CDB$, $\rm ABC$, and $\rm BAC$. 
This makes it easy to calculate coordinates.

$$
{\rm T}_5 = \left(\frac{-ar_1}{c}, \frac{-br_1}{c} \right),  \quad {\rm T}_6 =  \left( \frac{ar_3}{c},\frac{br_2}{c} \right)
$$

We are now in position for a spectacular result. The transformation 
$
\rho : (x, y)  \to  \  r_1 - \ y, \  \; r_1 - \ x \ 
$ is a self-inverse operation that fixes ${\rm T}_1, {\rm T}_2,$  ${\rm T}_3, {\rm T}_4$ 
and exchanges ${\rm T}_5$ with $ {\rm T}_6$! In other words, it is a reflection in the 
line $\lambda $. Most of the claim is nearly trivial. Showing that 
$\rho :{\rm T}_5 \longleftrightarrow {\rm T}_6$ amounts to verifying the following.
\begin{equation}
 \frac{r_1}{r_2} =  \frac {b}{a + c},  \quad \frac{r_1}{r_3}= \frac{a}{b + c}
\end{equation}
Substituting $b  =  r_1 + r_3,  \ a  +  c  =  r_4 +  r_2, \  a  =  r_1 +  r_2, \ b  +  c  =  r_4 +  r_3$ 
and simplifying reduces both to $ r_2  r_3 = r_1  r_4$, which is just \textbf{(7)}. 
The reflection $\rho $ must therefore take the four-circle system $\alpha $ to 
another system $\beta $. 
\begin{figure}[htbp]
\centerline{\includegraphics[width=2.5in]{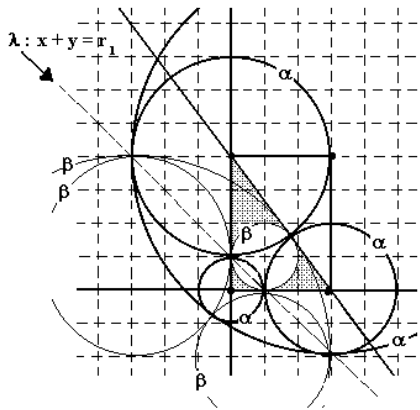}}
\caption{Reflection Duality of two tangent systems.}
\end{figure}
The result is pictured in \textbf{Fig.5} According to Coxeter {[5]}, the two 
sets of circles are Beecroft duals. It is no great effort to show that this 
combines the four triads into a single system $\alpha $ in a way 
that exhibits the companion orthogonal circles as the dual system $\beta $.

The four congruent triangles ${\text{ABC,  ABD,  ACD,  BCD}}$
are matched with the four $\beta $--circles; the matching circle is an in-circle or an 
ex-circle by \textbf{(1)}. Tangency for the in- and ex-circles of the same (scalene) triangle 
is not possible.

\begin {theorem}
The systems $\alpha$ and $\beta$ are Beecroft duals, congruent by reflection $ \rho : \ x + y  =  r_1,$ satisfying: \begin {enumerate} \item  Four of the tangency points lie on the line of reflection.  \item  The other two tangency points are symmetric to this line. \item  Systems $\alpha $ and $\beta$ share all six tangency points. \item  Each $\beta $-circle is orthogonal to three $\alpha $-circles. \end{enumerate}
\end{theorem}

There is something very satisfying about the eight interlocking circles 
illustrated in \textbf{Fig.5}. The complex of relationships here relies 
on the very elementary equations \textbf{(1)} and \textbf{(7)}. The alert 
reader may be curious at this point how those equations are related to the 
celebrated identity of Descartes (sometimes called the Soddy formula), which 
involves the sizes of four mutually tangent circles. To this identity we now 
turn.

\section{The Descartes Circle Theorem.} 
We rely heavily on [15]. Stated in essence in a letter (1643) by Descartes, 
this theorem was put more broadly and proved by Steiner (1826) and then others. Nobel Physicist 
F. Soddy made a delightful poem, and obtained a generalization to $n$ dimensions (1936).

First we define the \textit{oriented} curvatures $k_i$ which correspond to our radii $\cal R$.
$$ {\cal K }= [k_1,\ k_2, \ k_3, \ k_4] = [\tfrac{1}{r_1}, \ \tfrac{1}{r_2}, \ \tfrac{1}{r_3}, -{\tfrac{1}{r_4}}]$$

Basically the curvature is the inverse of the oriented radius. The oriented 
radius is either $r \ {\rm or}-r$,  where the unsigned radius is $r > 0$. In our context 
the negative sign is used exclusively for the fourth circle, and 
that is because it encloses the others.

The radial equations \textbf{(7)} can be rewritten as the following 
analogous equations \textbf{(11)}.
\begin{equation}
 k_2  k_3 +  k_1  k_4 =  0 , \quad   k_4 +  k_1 =  k_2 +  k_3 	
\end{equation}
The Descartes circle equation \textbf{(DCE)} is sometimes attributed to 
Soddy. It states that the sum of the squares of the oriented curvatures is 
half the square of the sum \textbf{(12)}.

\begin{equation}
[k_1^2 +  k_2^2 +  k_3^2 +  k_4^2] = \tfrac{1}{2}[k_1 +  k_2 +  k_3 +  k_4]^2
\end{equation}
This equation follows fairly easily (details omitted) from \textbf{(11)}.
The burden of the extra variables $k_i$ can be dispensed with. Since \textbf{(12)} 
is homogeneous, we can use \textbf{(8)} to replace $\cal K$ by $\textit{G}{\cal K} = [r_4,  r_3,  r_2,-r_1 ]$. 
Equally, the triple/triangle $[a, b, c]$ is rescaled to $\tfrac{1}{\textit{G}} [a, b, c]$ to replace 
$\cal K $ by $\textit{G}{\cal K}$. With given triple $\tfrac{1}{6} [3, 4, 5]$ we have 
$$
{\cal R} = \tfrac{1}{6} [1, 2, 3, 6], \quad  {\cal K} = [6,  3,  2,-1]
$$ 
This works. So does ${\cal K}  =  [20, 5, 12,  - 3]$ from $\tfrac{1}{6}[15, 8, 17]$:

$$
(36  +  9  +  4  +  1)  =  \tfrac{1}{2}(6  +  3  +  2  -  1)^2
$$
$$
(400  +  25  +  144  +  9)  =  \tfrac{1}{2}(20  +  5  +  12  +  3)^2
$$
Exercise: Find the four right triangles that are the sources for each of these 
solutions of the \textbf{DCE}.
$$
[\sqrt 2 + 2,  \sqrt 2, \sqrt 2, \sqrt 2 - 2 ] \quad  [ \sqrt 3 + 3, \sqrt 3 + 1, \sqrt 3 - 3,  \sqrt 3 - 1], 
$$
$$[\sqrt 5 + 3, \sqrt 5 + 1, \sqrt 5 - 1, \sqrt 5 - 3], \quad [ 15, 10, 3,  - 2 ]
$$
We will return to the topic of Descartes circles at the end. There is one 
other theme we need to pursue first.

\section {Pythagorean Triangles.}  
One might happen to notice that many of our examples are right triangles, all of whose sides 
are whole numbers, like $[3, 4, 5]$ or $[5, 12, 13]$.  A triple $[a, b, c]$
of positive integers satisfying \textbf{(6)} is called a {\textit{Pythagorean triple}}. 

Our equations \textbf{(1)}, \textbf{(2)} show that when the radii are integers, so 
are the sides of the (arbitrary) triangle. They also imply that when the sides $a, b, c$ 
are integers, \emph{and} at least one radius is an integer, then the other three radii are 
the sum/difference of two integers.

It remains to consider when the sides are integers, but the radii are 
half-integers, by \textbf{(2)}. For instance, the 
{1-1-1-}equilateral triangle has 
$$
{\cal R} =  [\tfrac{1}{2}, \ \tfrac{1}{2},\ \tfrac{1}{2}, \ \tfrac{3}{2}]
.$$ Doubling the sides \textit{and} the radii then gives an integer solution in which every 
radius is \textit{odd}.

For right triangles this cannot occur by \textbf{(7)}. Odd numbers are of two types: 
either $4k+1 $ or $4k-1 $. How many $r_i$ are of the second type? The first equation 
of \textbf{(7)} has a ``balance'' that requires an odd count, but the second equation requires 
an even count. Thus we conclude that any Pythagorean triple $[a, b, c]$ must correspond to 
integer radii $\cal R$ and conversely.

Besides, if the radii are multiplied by a factor $k > 1$, the same factor is applied 
to the sides, and \textit{vice versa}. This means that we can stick to \emph{primitive} 
triples, which have \textbf{GCD} equal to one. Clearly this is true for the triple if it is
 true for $\cal R$, and conversely.

However, our desire to employ $\cal R$ as the ideal set of parameters for a triangle is 
severely challenged in the right triangle case, since \textbf{(7)} must be satisfied! 
One strategy that can be used is to examine the restricting equations very closely. 
From \textbf{(7)} we remove the redundant $r_4$ and study $r_1  \cdot (r_1 +  r_2 +  r_3) =  r_2  \cdot r_3$. 
The first step is to multiply out and rearrange the terms: $r_2  \cdot r_3 - r_1  \cdot r_2  - r_1  \cdot r_3 -  {r_1}^2 =  0$. This very nearly factors (the last sign is wrong), but completing the product 
anyway gives us a useful equation.
$${2r_1}^2 =  (r_2 -  r_1) (r_2 -  r_1)$$ \begin{center}  {or} \end{center}\begin{equation}
	\tfrac{1}{2}(a  +  b  -  c)^2 =  (c  -  b)(c  -  a)
\end{equation}
The second equation is obtained by substitution from \textbf{(2)}. If we introduce 
variables $e, m, n $ for the three parentheses of the second equation, then \textbf{(14)} results.
\begin{equation}
a  =  m  +  e,  \quad    b  =  n  +  e, \quad  c  =  m  +  n  +  e,  \quad      e^2 =  2mn
\end{equation}

The essence of \textbf{(14)} was introduced a long time ago, by L. E. 
Dickson [6]. We agree with Gerstein [9] that it has been long overlooked; witness 
its rediscovery by Teigen and Hadwin [20] and others. We deal with this equation more 
fully elsewhere [3], showing that it can function as the source (almost effortlessly) of 
nearly all good and basic facts about Pythagorean triples. But perhaps the same is 
true of the circles approach in this article.

Klostergaard [13] and McCullough [17] use equivalent variables and equations to carry out 
enumerations. L\"{o}nnemo [16] uses it to construct the celebrated ``family tree'' 
(introduced in 1963 by Barning [2] or in 1970 by Hall [11]). Preau [18] uses a similar 
transformation to rediscover the tree. Alperin [1] also rediscovers it. We shall shortly 
give entirely new methods employing circles!

First we sketch a general algorithm for selecting the radii for a primitive 
triple. One example will suffice. Let $r_1 =  3$, say. Write $2{r_1}^2$ as a product of 
prime powers and one: $18  =  (1)(2)(9)$. With (1) and (2) in different camps, we may 
assign (9) to either camp; thus $18  =  (1)(18),$ or $18  =  (9)(2)$. The original `three' 
is then added to both factors, and the sum $r_4$ is computed. Thus $r_1 =  3$ leads to two possibilities:
$$
{\cal R} = [3, \  1 + 3, \  18   + 3, \ 28]  \to {\cal T} = [ \ 7, 24, 25]
$$
$$
{\cal R} = [3, \  9  + 3, \   \; 2   + 3, \ 20]  \to {\cal T} =[15, \ 8, 17]
$$
If we eliminate the restriction to primitive triples, it is easier to start by picking any value for
$$
m  =  c  -  b  =  r_1 -  r_2 .
$$ 
Then write $2m  =  f^2 g$, where $f^2$ is the largest possible square factor (and consequently 
$g$ is square-free). Finally, $n  =  g h^2 $, and $e  =  fgh$. This is a small simplification of 
McCullough's enumeration technique {[17]}. Did we mention that $e  =  2r_1$ is a diameter? 
(Yes, in the first section.)

We adopt the convention that when $[a, b, c]$ is a primitive triple, 
side $a$ is odd. Equivalently, the radii alternate between 
even and odd: we already know that at least one radius is even, but the 
second equation of \textbf{(7)} guarantees a second, and the first equation 
tells us two are odd (three even implies four even, contradicting 
primitivity). The second equation (again) makes one of $r_1, r_4$ even and one odd, 
likewise for $r_2, r_3$. We exchange $r_2, r_3$ if needed (also $a, b$) to get 
alternation between even/odd.

\section{Raising the Standard. } 
As further proof of the utility of circles, we now extract the age-old standard solutions 
for Pythagorean triples. The key is a ``minor'' lemma from \textbf{Section 3}, namely \textbf{(10)}. 
Here it is again, enlarged a bit, courtesy of the identities, \textbf{(6)}, \textbf{(7)}.
\begin {equation} 
	\begin{array}{ccc}
  	{\tfrac {r_1}{r_2} =   \tfrac {r_3}{r_4} =  \tfrac {b}{a  +  c} =  \tfrac {c  -  a}{b}}  \\ 
   	{} \\ 
   	{\tfrac {r_1}{r_3} =  \tfrac {r_2}{r_4} =  \tfrac {a}{b  +  c} =  \tfrac {c  -  b}{a}}  \\
\end{array}
\end{equation}

We avail ourselves of the opportunity to limit consideration to primitive 
triples. Hence the sequence $(r_i)$ has no non-trivial common 
factor, and alternates in parity. Consequently, the second group of equal 
fractions includes one with \textit{odd }numerator and \textit{odd }denominator.

We denote the reduced fraction in the first group by $q/p$, 
and the reduced fraction in the second group by $ q'/p'$. The pair $ (q,p)$ and the 
pair $ (q' \!,p')$ are both relatively prime, and the last is a pair of odd numbers.

It then works out that $r_1 =  xq, \;   r_2 =  xp, \; $ and $ r_3  = yq, \; r_4 = yp$, where 
$x$ is the \textbf{GCD} of $\ r_1  , r_2$ and $ y$ is the \textbf{GCD} of $r_3, r_4$. We note 
that any positive factor $k$ of $x,y$ is a factor of all four radii. Hence the pair $(x,y)$ 
is relatively prime. Finally $ \tfrac{r_1}{r_3} = \tfrac{r_2}{r_4} = \tfrac{x}{y} = \tfrac{q'}{p'}$ 
establishing $x  =  q', \;  y  =  p'$ and hence \textbf{(16)}.
\begin{equation}
r_1 = qq', \quad r_2 = pq', \quad r_3 = qp', \quad r_4 = pp'
\end{equation}
The key numbers $ q, \ p, \ q',\ p'$ are in fact pairwise relatively prime, 
since any factor $k$ that divides two of them must divide three radii, and 
by \textbf{(7)} all four. There is much more. First of all, a right triangle 
with legs $c + a, b$ is easily shown to have one acute angle equal to half an 
acute angle of right triangle $[a, b, c]$. Checking details shows that $q/p$ and $q'/p'$ 
are the half-angle tangents for $[a, b, c]$!

If we substitute from \textbf{(16)} into equation $r_1 +  r_2 +  r_3 =  r_4$ the result 
can be ``solved'' in either of the following two ways.

$$
\tfrac {q}{p} = \tfrac {p'  -   q'} { p'  +  q'} \ , \quad \tfrac {q'}{p'}  = \tfrac{p  -   q}{p  +  q}     
$$
Since we know that $ q',p'$ are odd, we can put $qz  =  \tfrac {1}{2}  (p'  -  q'),\; pz  =  \tfrac {1}{2}(p'  +  q')$. But then $ z(p  -  q)  =  q',\; z(p  +  q)  =  p'$, implying that $z  =  1$. And so
\begin{equation}
p - q  =  q', \quad  p + q  =  p', \quad \tfrac {1}{2} (p' - q')  =  q, \quad  \tfrac {1}{2} (p' + q')  =  p .	
\end{equation}
We can \textit{even} conclude that in the pair $(q,p)$ one is odd, and one is 
even. 

Using \textbf{(17)}, it is amusing to present the sequence ${\cal P} =  [q', q, p, p']$ as a 
Generalized Fibonacci sequence! In other words, $p  =  q' + q \ $ and $p'  =  q + p.$ 
One may pick, in fact, $q'$ to be any odd positive integer, and $q$ any 
positive integer which is relatively prime to $q'$. For example:
$$
{\cal P} \ =  [1,\ q,\ 1 + q, \ 1 + 2q]  \to {\cal R} \  =  [q,\; 1 + q,\; q(1 + 2q),\; (1 + q)(1 + 2q)]
$$
$$
 \to {\cal T} = [1 + 2q, \ 2q(1 + q), \ 1  +  2q(1 + q)]
$$
or with $q$ not divisible by three, 
$$
{\cal R} \ =  [3,\ q, \ 3 + q, \ 3 + 2q]  \to {\cal R} \  =  [3q,\ q(3 + q), \ 3(3 + 2q), \ (3 + q)(3 + 2q)]
$$
$$
 \to   {\cal T} = [3(3 + 2q), \ 2q(3 + q),\ 9 + 6q + 2q^2]
$$
More generally, we get, combining $[a, b, c]  =  [r_1 +  r_2, \; r_1 +  r_3, \; r_2 +  r_3]$ 
with \textbf{(16)}, \textbf{(17)}:
\begin{equation}
[a, b, c]  =  [p^2  -  q^2 , \; 2pq, \; p^2 +  q^2]
\end{equation}

\begin{equation}
[b,a,c] =  [\tfrac{1}{2}   (p'^2  +  q'^2 ), \  p'q', \  \tfrac{1}{2}   (p'^2  +  q'^2 )]
\end{equation}

\begin{equation}
[a, b, c]  =  [p'q',\;  2pq, \; pq'  +  qp']
\end{equation}

\begin{equation}
[a, b, c]  =  [p'q', \; 2pq, \; pp'  -  qq']
\end{equation}

Here \textbf{(18)}, \textbf{(19)} are the two standard solutions, the first in terms 
of $q,p$ alone, the latter (with $ a,b$ switched) in terms of $q',p'$ alone. We dare 
to add two mixed forms \textbf{(20)}, \textbf{(21)}, which seem to us to be the 
most generally useful! (Also, note the area is simply $\textit{G}  = qq'pp'$)!

Both standard solutions are found in Sierpinski {[19]}, in Dickson 
{[7]}, in Eckert {[8]}, and doubtless other places. One or the 
other, usually with an added factor of $k>0$ to allow for 
non-primitive solutions, are found in sources too numerous to list. The 
mixed forms are close in form and spirit to ``Generalized Fibonacci'' 
sequences, for example Koshy {[14]} or Horodam {[12]}, but are 
bit more symmetric than any that we have seen. The results of this section 
are recapitulated in the following theorem.

\begin{theorem}
\textbf{\textup{(}Radius Decomposition\textup{)}.}
Let  $[a, b, c]$ be a primitive triple with radii 
${\cal R} =  [r_1 , \ r_2 ,\  r_3,\  r_4.\ ]$ Then parameters $p,\ q,\ p', \ q'$ are determined such that, \textup{(16)}, \textup{(17)} are valid, and also \textup{(18)}, \textup{(19)}, \textup{(20)}, \textup{(21)}. Moreover, the fractions of \textup{(15)} are the half-angle tangents of the right triangle.
\end{theorem}
\section{Tree Rings (More Circles). }

Recall that ${\cal P}=[q', \ q,  \cdots  ]$ can be completed by ``Fibonacci addition'', 
$
{\cal P} =[q', \ q, \ p,\  p'],
$
 and we only require that $q'$ is odd, and $q',p$ are relatively prime. The four radii 
$
r_1 ,\ r_2 ,\ r_ , \ r_4
$
are found by multiplying either end value by either middle value-- this is exactly the 
result of expanding the product
$
(q'  +  p')(q  +  p)
$
Consider the possibility of modifying $\cal P$ by exchanging $q$ with $p$, or $q'$ with 
$p'$, or both. The result is then one of these.
$$
{\cal P}_1 =  [q', p, \cdots], \quad       {\cal P}_2 =  [p', q, \cdots], \quad    {\cal P}_3 =  [p', p, \cdots].
$$

The final two values are \textit{replaced} by computing with Fibonacci addition, and in each 
case, a valid sequence is obtained. 

We ask, what is the in-circle radius for 
${\cal P}_1,  \  {\cal P}_2, \  {\cal P}_3 $? It is found as the product of the first two members, namely $q'p, \ p'q, \ p'p.$ 
But these three values are just $r_2, r_3, r_4 $ (the ex-circle radii from the starting case $\cal P$)! 
In effect, we produce ${\cal P}_1,$ ${\cal P}_2,$ or ${\cal P}_3$ from $\cal P$ by ``promoting'' one ex-circle to the 
position of in-circle! Such \textit{promotion} always increases the size of the in-circle. Take the simplest example ${\cal P}=[1, 1, 2, 3]$ which has radii ${\cal R} =[1, 2, 3, 6].$ Promote $r_2$ to 
get ${\cal P}_1 = [1, 2, \cdots ]$, or promote $r_3$ to get ${\cal P}_2 = [3, 1, \cdots ]$, or 
promote $r_4$ to get ${\cal P}_3 =  [3, 2,\cdots ]$. Continuing with ${\cal P}_3$, we find ${\cal P}_3 = [3, 2, 5, 7],$ then 
multiply out $(3  +  7)(2  +  5)  =  (6  +  15  +  14  +  35)$.  Finally, add $(6 + 15),$ $(6 + 14),$ 
$(15 + 14)$ to get the triple $[21, 20, 29]$.

The cancellation of a ``promotion'' is a ``demotion''. We claim that 
any $\cal P$ other than ${\cal P} = [1, 1, 2, 3]$ has a unique demotion! 
Demotion \textit{decreases} the in-radius.

\begin{theorem}\textbf{\textup{(}Demotion\textup{)}.} Let 
$
{\cal P}=[x, y, \cdots ]
$
be a valid sequence other than $[1,1,2,3]$. Then there are three uniquely determined Fibonacci-rule sequences of form $[x, *, y, * \ ], \ [*, y, *, x \ ], \ \ [*, *, y, x \ ].$
Exactly one of these has all entries positive. This is the demotion ${\cal P}^-$ of $ \ {\cal P}.$ Thus 
$\cal P$ can be demoted in exactly one way.
\end{theorem}

The proof of this theorem is so direct, we need only sketch the demotion 
algorithm and leave the rest to the reader. If ${x  =  y }$ then 
${\cal P} =  [x, x, 2x, 3x]$  $=  [1, 1, 2, 3]$, so assume $x  <  y$ or $x  >  y.$

If $x  <  y$ then $[x, *, y, *]  =  [x,  y - x,  y,  2y - x]$ works and the 
other two do not. If $x  >  y$ take $[*, *, y, x]$ and ``back up'' to $ [*, z, y, x].$ Note $z=y$ is 
excluded, since $x=z+y$ is odd! Then, if $z>y$, reorder to $[*,y,z,x]$ (sort into ascending order). 
Finish backing up. The result is then either $[*, *,  y,  x]   =   [2y - x,  x - y,  x,  y]$ or 
$[*,  y, *,  x]  =  [x - 2y,  y,  x - y,  x].$ Note the initial element is always odd.

Unlimited demotion (infinite descent) is impossible, for the in-radius is 
always decreasing. Therefore, any triple can be obtained in one way by a 
sequence of promotions, starting with $ {\cal T}  =  [3, 4, 5]$
, that is, ${\cal P}  =  [ 1, 1, 2, 3].$  The set of primitive triples then has the structure of an 
infinite ternary tree. At the top is the ``root'' ${\cal T}_0 =  [3, 4, 5].$  Each ``node'' 
is a primitive triple with three immediate successors. These successors are (A) left, 
(B) middle, or (C) right respectively as the promoted circle has size $r_3 , r_4 , r_2 .$  
For the root, this means $[15, 8, 17]$ is left, $[21, 20, 29]$ is middle, and $[5, 12, 13]$ is right. 
The repertoire of simple connection equations in \textbf{Theorem 2} makes it easy to go from one to another. 
For illustration, we compute $\cal P$-sequence $[3,2,5,7]$ and $\cal R$-sequence $[6,15,14,35]$ for the 
middle triple $ {\cal T}  =  [21, 20, 29]$.  Readers should produce enough examples to demonstrate 
this for themselves. 

\begin{cor}
\textbf{\textup{(}Barning/Hall Family Tree\textup{)}}. The set of primitive Pythag-orean triples can be 
seen as the nodes of a ternary tree. Demotion moves up one level to the ``parent'' and 
promotion moves down one level, producing three immediate successors ``children''. Only one triple $[3,4,5]$ has no parent.
\end{cor}

A certain amount of tedious computation is all that stands in the way of 
expressing promotion/demotion as a transformation of triples $[a,b,c]$, or a transformation of half-angle tangents $(q,p)$ , $(q',p')$. See for example [3]. Here the following will suffice.

\begin{cor}\textbf{\textup{(}Half-angle tangents\textup{)}}. 
Let $ \tfrac{x}{y}$ be either half-angle tangent for $[a,b,c],$ but not $\tfrac{1}{2}$ or $\tfrac{1}{3}$. The fraction is demoted to $\pm (\tfrac{x}{y - 2x})^\pm .$ The first $\pm 
$ is chosen to avoid negative quantities; the second $\pm $ to avoid fractions greater than one. {Examples}: 
$$\tfrac{3}{10} \rightarrow \tfrac{3}{4},  \quad \tfrac{3}{8} \rightarrow \tfrac{3}{2} \rightarrow \tfrac{2}{3}, \quad \tfrac{5}{6} \rightarrow \tfrac{5}{2} \rightarrow \tfrac{5}{4} \rightarrow \tfrac{4}{5}$$
\end{cor}

\begin{cor}
\textbf{\textup{(}Reversed Barning-Hall Transformation\textup{)}}. If triple $[a,b,c]$ is not equal to $[3,4,5],$ it may be demoted to
$$
 |a - 2e|, \ |b - 2e|, \ c - 2e  =  | -a - 2b + 2c|, \ | -2a - b + 2c|, \ -2a - 2b + 3c
$$
$Example: \ \left[33, 56, 65\right] \to \left[\ |33-48|, \ 56-48, \ 65-48\right]  =  \left[|-15|,\  8, \ 17\right].$
\end{cor}

Comments: See \textbf{(14)} for $e = a + b - c$. As L\"{o}nnemo 
{[16]} points out,  \textbf{(14)} works just as well with $e$ changed to $-e$. This clue leads to the 
transformation of \textbf{Corollary 2} except for the absolute value operations.

The reader may wish to work out the corresponding promotions implied by these corollaries: 
ratio $\tfrac{q}{p}$ promotes to $\tfrac{q}{p+2q}, \ \tfrac{p}{2p+q}, \ \tfrac{p}{p-2q}$; 
triple $[a,b,c]$ promotes to
\begin{equation}
[  \pm  a   \pm  2b  +  2c, \quad \pm 2a  \pm  b  +  2c, \quad \pm 2a  \pm 2b  +  3c ]
\end{equation}
where in all three terms the pattern is the same one: $(-++)$ or $(+++)$ or $(+-+)$.

If the triple is taken as a column vector $\rm X$, then 
promotion is just a matrix linear transformation: $\rm X \rightarrow AX$ or 
$\rm BX$ or $\rm CX$. Here $\rm B$ is 
square with $(1,1,3)$ on the diagonal and 2 
elsewhere, and changing the signs in the first (the second) column gives 
$\rm A $ (gives $\rm C$).

We have already cited the independent discovery of the ternary tree, seven 
years apart, by Barning and Hall. Eckert {[8]} alludes to the ideas 
of \textbf{Corollary 2}. See Dickson {[7]} for the fact that Fermat 
essentially found transformation ${\cal P} \rightarrow {\cal P}_3$, but (apparently) not 
either of the other two. Of all the different ways to construct the tree, we 
believe our circle method is the simplest, and the most transparent.

We have yet one more original way to construct the tree via circles! The key 
is a beautiful circle from classic geometry, studied by Feuerbach and many 
others. We think of it in this context as the ``family circle.''

In preparation, we recall \textbf{Fig.4} and \textbf{Fig.5}. The four 
tangent circles in cluster $\beta$ are ``attached'' as 
equi-circles to congruent right triangles $\rm ABC,$ $\rm ABD,$ $\rm CDA,$ $\rm CDB.$ Focus on $\rm ABC$. The equi-circles for $\rm ABC$ are ${\cal C}(\rm I_i,  r_i)$ where 
$$
\rm I_1 =  (r_1,  r_1), \quad \rm I_2 = (r_2,  -r_2), \quad \rm I_3  = (-r_3,  r_3), \quad \rm I_4 = (r_4,  r_4).
$$
The four symmetries of rectangle $\rm ADBC$ are obtained from the following:
$$
{
\begin{array}{rlrl}
\sigma : x & \to \     x,	 &  \ \tau : x &  \to \ a - x \\ 
          y & \to \  b - y,  &           y & \to \ y 
\end{array}
}
$$
These are reflections, vertical and horizontal. Recall also $$ \rho : (x, y)  \to (r_1 - y,\ r_1 - x).$$ 
The four $\alpha$-circles ${\cal K} =  {\cal C}(\rm C, r_1),  \ {\cal C}(\rm B, r_2), \ {\cal C}(\rm A, r_3), \ {\cal C}(\rm D, r_4 )$ are mapped to the above equi-circles by 
${\cal K}  \to \rho {\cal K}, \   \tau \rho {\cal K}, \   \sigma \rho {\cal K}, \ \sigma \tau \rho {\cal K}.$

The center $\rm M  = (\tfrac {a}{2}, \ \tfrac {b}{2})$ of the rectangle is also the midpoint of 
diagonals $\rm \overline{AB}$ and $\rm \overline{CD}$. Circle ${\cal C}(\rm M, \tfrac {\rm c}{2})$ is seen 
to be the circum-circle of the rectangle, and also of the four right triangles. Let 
$\rm F = \tfrac {\rm M}{2} =  (\tfrac {a}{4},\tfrac {b}{4}).$  Circle ${\cal C}(\rm F, \tfrac {c}{4})$ 
has $\rm \overline{CM}$ (one radius of the circum-circle) for a diameter. It is also a circum-circle 
for a ``shrunken'' rectangle $\rm A'B'C'D'$:  We note that $\rm A'B'D'$ are the midpoints of the 
sides of $\rm ABC$. Hence this last circle with center $\rm F$ is the {\textit{nine-point circle}  $\cal N$ 
studied by Feuerbach and others.

In \textbf{Fig.6} the equi-circles, the circum-circle, and the nine-point circle are shown. 
Note point $\rm L$, the foot of the altitude from $\rm C$. The $\overline {\rm AB}$ diagonal line, 
and the altitude $\overline {\rm CL}$ have equations $bx + ay = ab, \ ax = by .$ Hence the foot 
of the altitude is $$\rm L = ( \tfrac {ab^2}{c^2 }, \tfrac {a^2 b}{c^2}).$$
	\begin{figure}[htbp]
	\centerline{\includegraphics[width=3.2in]{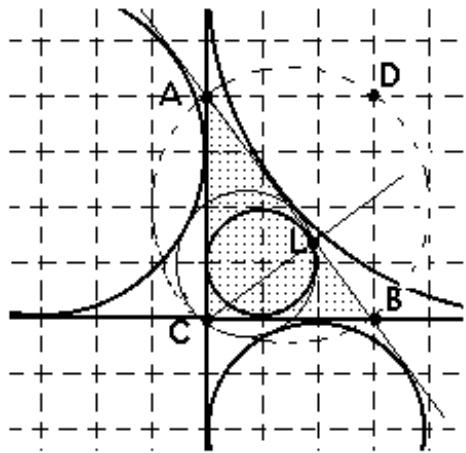}}
		\caption{In-circle, ex-circles: all tangent to nine-point circle $\cal N.$}
	\end{figure}
Origin $\rm C$ is the foot of the other two altitudes, and also a point of tangency between the circum-circle and its ``half'' brother $\cal N$. As diameter $\rm \overline{CM}$ of $\cal N$ is also hypotenuse for right triangle $\rm CLM$, the point $\rm L$ is on circle $\rm M$. It is an interesting exercise to verify this by finding the side lengths, and using the Pythagorean theorem. A helpful datum for this pursuit is the Pythagorean right triangle with hypotenuse $c^3$ and legs equal to the absolute value of the real and imaginary parts of the third power $(a + bi)^3$.  In \textbf{Fig.6 } circle $\cal N$ seems tangent to the four equi-circles -- no surprise to geometers! We shall verify this beautiful result of Feuerbach, and in the process, make a wonderful and serendipitous discovery! 
The key is to compute the (vector) difference of the centers of two circles, 
and the sum/difference of the radii. Should the values (aside from sign) 
turn out to be the legs and hypotenuse of a right triangle, the circles are 
tangent. Thus $$ \rm I_4 - F = (r_4 - \tfrac {a}{4} , \ r_4 - \tfrac {b}{4}) = \tfrac {1}{4}(a+2b+2c, \ 2a+b+2c)$$
$$\text{Sum of radii:} \ r_4 + \tfrac {c}{4}= \tfrac {1}{4}(2a+2b+3c).$$
 The ``fingerprints'' of \textbf{(22)} have shown up (middle case) --- the segment 
$\rm \overline{I_4F}$ is the hypotenuse of a quarter-sized copy of the middle child on the tree! 
By similar reasoning, the segments $\rm \overline{I_2F},\ \rm \overline{I_3F}$ go with the other two 
children (quarter-sized).

\textbf{Fig.7} magnificently illustrates the production of the 
three ``infants''. But it does nothing to explain why the ex-circle used to 
construct one of the ``infant'' triangles is exactly four times the size of 
the in-circle of that triangle, a fact we know because of the study of 
promotion at the beginning of this section.

\begin{figure}[htbp]
\centerline{\includegraphics[width=3.2in]{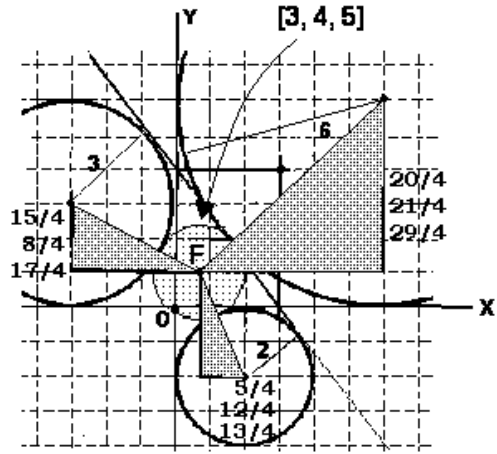}}
\caption{Children of $[3,4,5]$ emerge quarter-sized from tangencies.}
\end{figure}

But what about the other two circles tangent to circle $\cal N$, namely the in-circle with 
center $\rm I$ and the circum-circle with center $\rm M$? Segment $\rm \overline{FM}$ is 
the hypotenuse of a quarter-sized copy of the main triangle, and $\rm \overline{FI}$ is the 
hypotenuse of the ``parent'' of the main triangle (all quarter-sized, except, of course, the 
root triple $[3,4,5]$, which has a degenerate parent $[0,1,1]$. Truly the nine-point circle of 
a Pythagorean triangle is a ``family'' circle!

One us (Price) was studying the results of \textbf{Fig.7} and found a 
convenient way to produce the immediate family (one predecessor and three 
successors) of a given triple ${\cal T} = [a,b,c]$. Compute $t = 2(a+b+c)$ and then check 
that $[t-a, \ t-b, \ t+c]$ is a positive triple with a larger hypotenuse, with the same 
\textbf{GCD} (so both are primitive if either is). Note that $t$ is double the perimeter of triangle $\cal T$.

Now suppose that exactly one of $a,b,c$ is negative. The three 
cases may be distinguished by the sequence of signs: $(-++), (+-+), (++-).$  We add $(+++)$ 
for the all-positive example of the previous paragraph. For $(-++), (+-+),$ 
$t/2$ is not the perimeter of course, but the result is still a triple with larger hypotenuse. 
In the last case $(++-)$ we must take absolute values to get a positive triple: $ |t-a|, |t-b|, |t+c|$. 
This time, the hypotenuse is smaller. 

Changing two or three signs of $a,b,c$ only duplicates the four 
cases already seen, plus some added sign changes. We can now recognize the 
left, middle, and right child of a triple as the same as the outcome of 
cases $(-++), (+++), (+-+)$, and also recognize the fourth case (++-) as the 
predecessor! In hindsight, the method is close to the generalized Fermat 
method of Eckert [8].

Finally we are ready to wrap up the fireworks with a finale, one with lots 
and lots of circles!

\section{A Descartes-Pythagoras Reunion. }
In \textbf{Section 4} we derived 
the Descartes Circle Equation for a very special case. With [15] as 
reference, we now claim it for any set of four circles, tangent at six 
distinct points.

Our immediate goal is to examine certain integer solutions. We already have 
the {\textit{basic solution}} ${\cal K}' = [r_4, r_3, r_2,-r_1]$ belonging to 
rescaled triangle $\tfrac{1}{\textit{G}}[a,b,c]$. Provided that $[a,b,c]$ is 
a primitive Pythagorean triplet, integer quadruple $\cal K'$ is 
also primitve.

The triad of tangent circles with curvatures $r_4, r_3, r_2$ are tangent 
to the surrounding circle with curvature $-r_1$, but are also 
all tangent to a different circle placed in the midst of them. Let that 
smallest circle have radius $\textit{G}/r_0$ (or radius $r_0$ for $[a,b,c]$). The 
curvature (positive) turns out to be $\textit{G}/r_0 = \textit{G}k_0 = 4r_4 - r_1$. For the 
triple $[3,4,5]$ we find $\textit{G}k_0 = 4(6) - 1 = 23$. It is easy to check 
that $[6,3,2,23]$ satisfies the \textbf{DCE}. 

Let $f,g,h$ be any three curvatures, $f,g$ positive, and 
$\tfrac{1}{f} + \tfrac{1}{g} + \tfrac{1}{h} > 0.$ A corresponding triad of three 
circles tangent at distinct points may be constructed; the third circle encloses the others 
in case $h$ is negative. We may substitute $f,g,h$ into \textbf{(12)} \ to obtain a quadratic, $
x^2 - {\rm S}x + {\rm T = 0},$ with two real solutions $x_1, x_2$, such that ${\rm S} = x_1 + x_2$ and $ {T} = x_1 \cdot x_2$.

\begin{figure}[htbp]
\centerline{\includegraphics[width=2.3in]{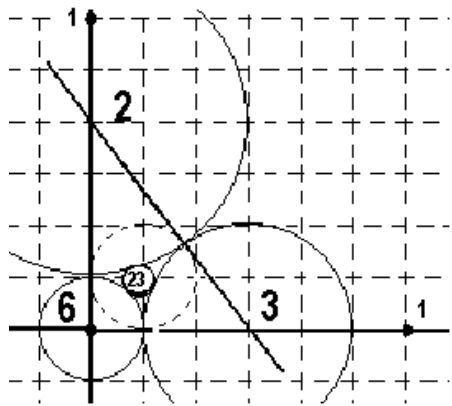}}
\caption{A different ``in'' circle for the triangle $\tfrac{1}{6}[3,4,5]$.}
\end{figure}

It is both clear and commonly known that two different fourth circles exist, 
tangent to the given three. Evidently $x_1,x_2$ are these two curvatures. Exchanging 
one circle for the other is called \textit{Descartes reflection}. We also find
\begin{equation}
x_1 + x_2 = {\rm S} = 2(f + g + h).
\end{equation}

Thus $[6, 3, 2,-1],$ $[6,3,2,23]$ satisfy \textbf{(23)}, since 
$23 + (-1) = 2(6 + 3 + 2)$. We see from \textbf{(23)} that Descartes reflection can be 
applied to any one circle out of a set of four tangent circles. Moreover, if all four 
circles have integer curvature, so does the new circle produced by reflection! 

Descartes reflection is a special case of circular inversion. Consider the 
in-circle (dashed) in \textbf{Fig.8}. Reflection in this circle must fix 
the three circles to which it is orthogonal $(2,3,6)$, and swap 
the other two $(23, -1)$ because tangencies are preserved. Much 
more information is available in [15] or [10], including 
simple means of locating the centers of reflected circles. 
Our interest in Descartes reflection is more limited in this paper. 

Suppose we start with \textbf{Fig.8} rotated $90^ \circ$ clockwise. Then 
by line reflections, or by $\rm D$-reflections of individual circles, we get 
\textbf{Fig.9(a)}, an arrangement of nine circles.

\begin{figure}[htbp]
\centerline{\includegraphics[width=3.9in]{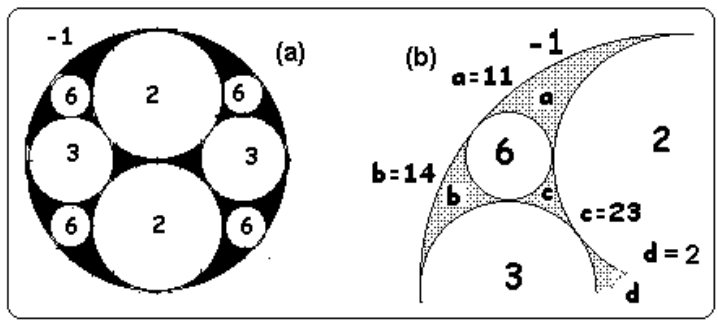}}
\caption{Descartes quadruple $[6 , 3, 2,-1]$ leads to a packing of circles.}
\end{figure}

\textbf{Fig.9(b)} is an enlarged view of the upper left of \textbf{9(a)}. 
When circles fill the gaps $\rm (a,b,c,d)$; the curvatures $11,14,23,2$ are found by using \textbf{(23)}.
$$
\begin{array}{rl}
	3+x = 2(6 + 2 - 1), \quad & 2 + x = 2(6 + 3 - 1), \\
	-1 + x = 2(6 + 3 + 2), \quad & 6 + x = 2(3 + 2 - 1). 
\end{array}
$$
Since all four calculations start from $[6, 3, 2,-1]$, it is possible 
to package them together. Replicate twice the sum $[20, 20, 20, 20]$, and 
also multiply each curvature by three, to get $[18, 9, 6,-3]$. Subtract as vectors: 
$$[20-18, 20-9,\ 20-6,\ 20+3] = [2, 11, 14, 23].$$ Each of these values can 
replace the original value in the same position.
$$
[ \underline{2},\ 3,\ 2, -1], \quad [6, \underline{11},\ 2, -1], \quad [6,\ 3, \quad \underline{14},  -1], \quad [6,\ 3, \ 2, \ \underline{23}]
$$
Clearly this process can be continued indefinitely, filling the interior of the 
``negative'' circle with an infinite ``foam'' of circles of increasing 
integer curvature. This is an \textit{integral Apollonian packing} (\textbf{IAP}). In both [10] and 
[15] is a beautiful picture which carries our \textbf{Fig.9} to many more levels. We will shortly present 
a detailed picture of our own.

The authors of [10], [15] offer no connection between integer 
packings and integer right triangles. But in \textbf{Fig.10(a)} see how 
the centers of the largest circles form four rectangles. The shaded right 
triangle, similar to $[3,4,5]$, is for us a nucleus for the 
packing.

\begin{figure}[htbp]
\centerline{\includegraphics[width=2.73in,height=1.25in]{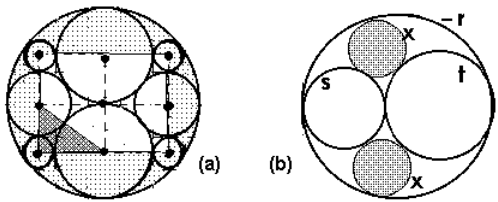}}
\caption{Packings from Pythagorean rectangle, or symmetric starting point.}
\end{figure}

There is little that is privileged about the $\tfrac{1}{6}[3,4,5]$ triangle. 
We can get different Apollonian packings from any rescaled primitive triple 
$\tfrac{2}{ab}[a,b,c]$. It will always be an \textbf{IAP} with \textit{two} rectangles, and 
bilateral symmetry. The children of $[3,4,5]$, rescaled, are
$$
\tfrac{1}{60}[15,8,17], \quad \tfrac{1}{210}[21,20,29], \quad \tfrac{1}{30}[5,12,13].
$$

This gives us three more \textbf{DCE} solutions ${\cal K}'$ in which we now {\textit{sort}} the 
entries to put it in the form $[-a, b, c, d]$ with $0 < a < b \leq c \leq d$:

$$
[-3, 5, 12, 20], \quad [-6, 14, 15, 35], \quad [-2, 3, 10, 15].
$$

The problem is that different sequences ${\cal K}'$ may belong to 
the same \textbf{IAP}. Given two sequences, each determines a unique 
\textbf{IAP}, but how can one tell if the packings are the same or distinct? 
A good question, but first we show how to generate many more packings! That 
will bring us near to the end of our study.

As we prepare to ``pack it in'' we return ``full circle'' to a figure 
reminiscent of \textbf{Fig.1}, namely \textbf{Fig.10(b)}. 
Suppose that $r, s, t$ are positive integers, and that three 
tangent circles with oriented curvatures $r, s,-t$ have 
their centers on a line. Since the radii form a sum, the (absolute) 
curvatures must make a harmonic sum: $$\tfrac{1}{t} = \tfrac{1}{r} + \tfrac{1}{s}.$$
We desire that the Descartes equation applied to ${\cal K}' = [-t, r, s, x]$ will 
have two identical integer roots $x: \ x = x_1 = x_2$. In that way we can 
generate an \textbf{IAP} with bilateral symmetry. The harmonic sum can be written in the form $(r-t)(s-t) = t^2$, which 
resembles \textbf{(14)}, and has this Diophantine solution (easily obtained).
$$
r = km(m+n), \quad r = kn(m+n), \quad t = kmn.
$$

Here $k$ is the \textbf{GCD} for $r, s, t$, assuming that $m,n$ are relatively prime. 
We want primitive solutions, so we put $k = 1$. To avoid duplication by symmetry, we put $m < n$
(let $m = n = 1$ be designated as \textit{isolated}). Again, serendipity! The condition for 
identical roots $x = u$ in the Descartes equation is automatically fulfilled, so that \textbf{(23)} becomes 
$u = r + s -t$, thence \textbf{(24)}.
\begin{equation}
	r = m(m+n), \quad s = n(m+n), \quad t = mn, \quad u = m^2 + mn + n^2
\end{equation}

\textbf{Fig.9(a)} clearly shows that the isolated solution ${\cal K} = [-1, 2, 2, 3]$ and 
solution ${\cal K}' = [-1, 2, 3, 6]$ are in the same packing ({6} reflects to {2} and vice versa). 
We therefore consider $m < n$. Suppose we have arbitrary primitive triple 
$[a,b,c] = r_1 + r_2, \ r_1 + r_3, \ r_2 + r_3$. Then generator pairs $(q, p)$ and $(q',p')$ can 
both take turns as $m,n$ in \textbf{(24)}, in which case the resulting Descartes solutions can be 
written as \textbf{(25)}.

\begin{equation}
	{\cal K} = [ \tfrac{-b}{2}, r_3, r_4, \tfrac{c+b}{2}], \quad {\cal K}' = [-a, 2r_2, 2r_4, 2c+a]
\end{equation}
Take triple $[7, 24, 25]$ for instance. From the generators $\tfrac{m}{n} = \tfrac{3}{4}, \tfrac{1}{7} $ we 
derive by \textbf{(24)} or \textbf{(25)} two bilateral solutions $[12, 21, 28, 37]$ and $[-7, 8, 56, 57^{==}]$.  Here $k^{==}$ indicates that a curvature $D$\textsl{-reflects} to itself. 
These agree with \textbf{(25)}. We also have the more natural $ {\cal K}' = [-3, 4, 21, 28]$.

There is another type of bilateral symmetry to consider. It splits via parity into two subcases, 
illustrated by the \textbf{DCE} solutions $[-4, 8, 9, 9]$ and $[-3, 5, 8, 8]$. A fairly direct analysis, similar to the above yields the following 
algebraic solution for $[-a, b, c, c].$
$$
a = 2k(t-k), \quad b = 2k(t=k), \quad c = t^2
$$
In this solution, we suppose that $\tfrac{t}{k} > 2$ is a reduced fraction. If $t$ is even, 
then divide $a,b,c$ by two. 

If the packing is identified by the simplest quadruple it contains (largest 
possible circles), called the \textit{root}, then quadruple $[-3, 4, 21, 28]$ given 
above is \textit{not} a root. Producing the root for this case is easy. Just reflect the 
largest number three times.
$$
[ - 3, 4, 21, 28] \ \to  \ [ - 3, 4, 16, 21] \  \to \ [ - 3, 4, 13, 16] \  \to \ [ - 3, 4, 12, 13]
$$

More generally, let ${\cal K} = [-u, r, s, t]$. Let $-u + r + s = v.$ Then if $t$ exceeds $v$, $\rm D$-reflection replaces $t$ with a smaller integer. Repeat until the the process ends, either with $[-a, b, c, d^{==}]$ in which $d$ reflects to $d$, or with $[-a, b, c, d^{><e}]$ where $d < e$ and $d$ 
reflects to $e$.

In [10] it is shown that this \textit{reduction} always leads  to a root quadruple 
$\cal K$ defining the \textbf{IAP}. Our process finds many -- see 
the table at the end of the article. One of the simplest packings (see \textbf{Fig.11}) comes from 
the root quadruple $[-2, 3, 6, 7^{==}]$. Isn't ``she'' lovely!

\begin{figure}[htbp]
\centerline{\includegraphics[width=3.10in,height=2.74in]{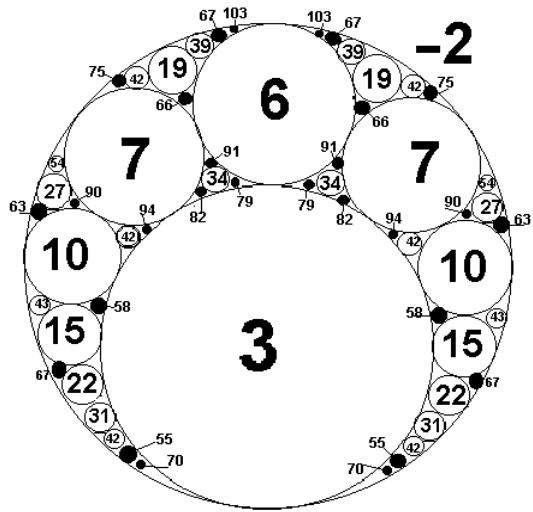}}
\caption{Apollonian packing $[-2,3,6,7]$ contains ``rectangle'' $[-2,3,10, 15].$}
\end{figure}

\section{Reverse Trends.} 
We have seen how to get many circle 
configurations out of a \textit{single} Pythag-orean triple. But now we glance at 
\textbf{Fig.12} to discover that one tangent circle array can be used to make more triples!
This is a simplified version of \textbf{Fig.8} showing the centers $\rm C, B, A, D, Z$ of five circles, each marked with a curvature $(1, 2, 3, 6, 23)$. 

\begin{figure}[htbp]
\centerline{\includegraphics[width=2.8in]{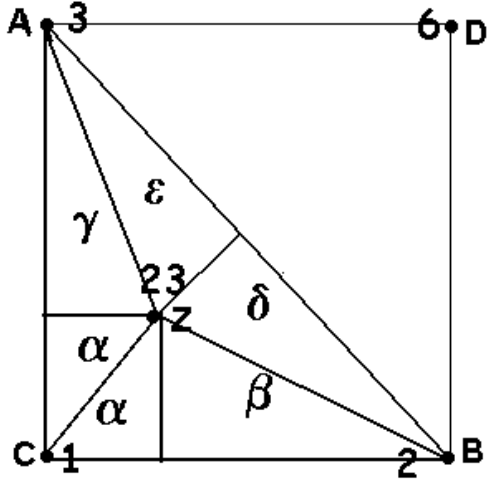}}
\caption{Dissection into rational right triangles.}
\end{figure}
By various techniques, one can (a) find the center $\rm Z$ of the (23)-circle, (b) join $\rm Z$ as shown to the vertices $\rm ABC$, and (c) drop perpendiculars to the sides $\rm \overline{AB},\ \overline{BC}, \ \overline{AC} $. This divides the triangle $\rm ABC$ into six smaller right triangles, two of them congruent $(\alpha)$.
In another article [4] we plan to analyze the general case which this diagram points to, and show that all the small triangles are rational right triangles. Specifically, in this example, the triangles 
$\alpha ,\beta ,\gamma ,\delta ,\varepsilon $ are similar to right triangles
$$
[21, 20, 29],  \quad  [5, 12, 13],  \quad [7, 24, 25],  \quad  [33, 56, 65], \quad   [117, 44, 125].
$$ 

The enterprising reader is invited to find appropriate multiples of these triangles and fit 
them together like a jigsaw to make a multiple of $[3,4,5]$.

\section{Table of integer DCE solutions}
\begin{figure}
\centerline{\includegraphics[width=5.18in,height=7.41in]{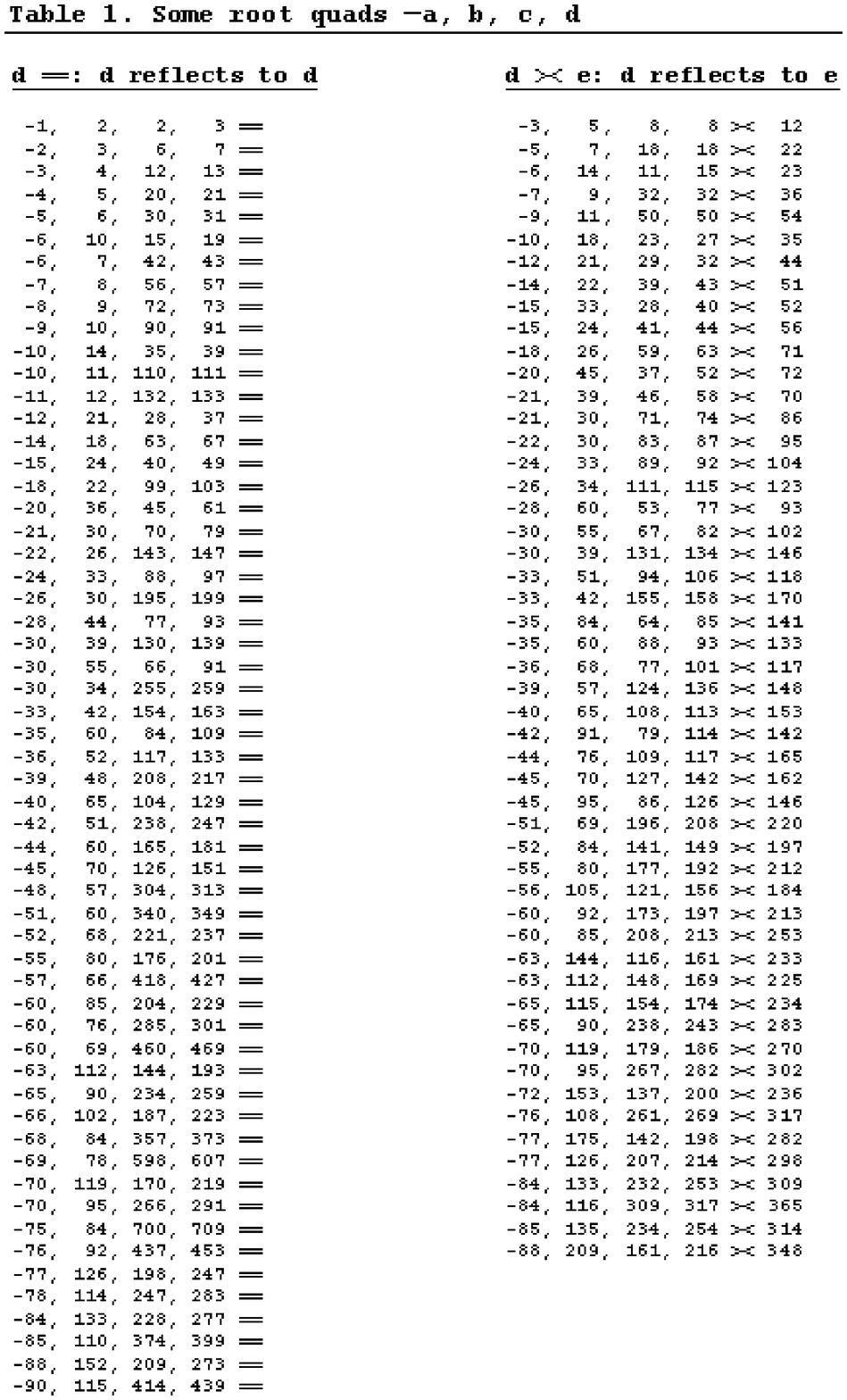}}
\end{figure}
The table contains only sorted sequences of form $[-a, b, c, d^{==}]$ in which $d$ is its own 
reflection  (i.e.$-a + b, + c = d$) or of form $[-a, b, c, d^{><e}]$ in which $d < e$ and $d$ reflects 
to $e$. The first category always gives \textbf{IAP}s that are bilaterally symmetric. All entries 
are ``root quadruples'' as in [10]. All are all obtained by reducing quadruples obtained from \textbf{PPT}'s. 

All solutions from \textbf{(24)} are recognizable in that $-a+b = d-c = m^2$ and $\tfrac{b}{c} = \tfrac{m}{n}$. This includes all in the first category. When one inspects a diagram such as \textbf{Fig.11}, the rectangle characteristic of a \textbf{PT} is spotted by placing horizontal lines through the centers of the (-2) and (3) circles, joining 
also the (10)s and the (15)s: $[-2, 3, 10, 15]$. Equivalent numerical tests are scarce. Certain patterns compel attention. We see the symmetric $[-m, m+2, n, n^{>< \ n+4}]$ already discussed. Similarly, we see the pattern $[-m, m+8, n, n+4]$ and solve it:
$$
m = 2(2k-1), \quad n = 2k^2 + 2k-1 .
$$
One more: $[- m,\ m + 9,\ n,\ n + 3]$ is solved by $m  =  3k,\ n  =  k2  +  3k  +  1.$ We have 
constructed over a dozen such one-parameter families. An added bonus is that these formulas have \textit{bona fide} instances not in the table! Clearly a more encompassing theory is to be sought! $${ }$$

\newpage

\textbf{Author Contact Information}
\begin{itemize}
\item \small Frank R. Bernhart, (\textit{bernhart@math.uiuc.edu})\newline
Math Dept Visitor, University of Illinois, Urbana, IL 61801.

\item \small H. Lee Price, (\textit{tanutuva@rochester.rr.com})\newline
83 Wheatstone Circle, Fairport, NY 14450-1138.

\end{itemize}

\begin{thebibliography}{\textbf{References}}

\bibitem{}  Roger C. Alperin, The Modular Tree of Pythagoras, Amer. Math. Monthly (112) (2005), 807-816.

\bibitem{}  F. J. M. Barning, On Pythagorean and quasi-Pythagorean triangles and a generation process with the help of unimodular matrices. (Dutch) Math. Centrum Amsterdam Afd. 
Zuivere Wisk. ZW-001 (1963).

\bibitem{}  Bernhart and H. L. Price, On Pythagorean Triples I, (in preparation)

\bibitem{}  Bernhart and H. L. Price, Pythagorean Triples by Dissection, (in preparation)

\bibitem{}  S. M. Coxeter, ``Introduction to Geometry'', John Wiley {\&} Sons, NY (1969, 2$^{nd}$ Ed. 1989)

\bibitem{}  L. E. Dickson, Lowest Integers Representing Sides of a Right Triangle, Amer. Math. Monthly {1} (1894), 6-11.

\bibitem{}  L. E. Dickson, ``History of the Theory of Numbers, Vol. II: Diophantine Analysis'', Chelsea Publ. Co., NY.

\bibitem{}  E. J. Eckert, Primitive Pythagorean Triples, College Math. Journal 23,5 (1992) 413-417.

\bibitem{}  Larry J. Gerstein, Pythagorean Triples and Inner Products, Math. Magazine 78 (2005), 205-213.

\bibitem{}  R.L. Graham, J.C.Lagarias, C.L. Mallows, A.R. Wilks, C.H. Yan, Apollonian Packings : Number Theory, J. Number Theory, 100, (2003) 1-45.

\bibitem{}  Hall, Genealogy of Pythagorean Triads, Math. Gazette 54, No. 390 (1970), 377-379.

\bibitem{}  F. Horodam, Fibonacci Number Triples, Amer. Math. Monthly 68(1961), 751-753.

\bibitem{}  Henry Klostergaard , Tabulating All Pythagorean Triples, Math. Magazine (1978), 226-227.

\bibitem{}  Thomas Koshy, Generalized Fibonacci Pythagorean Triples, Notes 86.61, Math. Gazette 86 (2002), 459.

\bibitem{}  J.C. Lagarias, C.L. Mallows, and A.R. Wilks, Beyond the Descartes Circle Theorem, Amer. Math. Monthly 109 (2002), 338-361.

\bibitem{}  H. A. L\"{o}nnemo, The Trinary Tree(s) underlying Primitive Pythagorean Triples, June 8, 2000. In: Cut the Knot, Interactive Mathematics Miscellany and Puzzles.A.Bogomolny(Ed.), http://www.cut-the-knot.org/pythagoras/ PT matrix.shtml

\bibitem{}  Darryl McCullough, Height and excess of Pythagorean triples, Math. Mag. 78(2005), 26-44.

\bibitem{}  Paul Pr\'{e}au, Un graphe ternaire associ\'{e} \`{a} l'\'{e}quation $\rm X^2 + Y^2 = Z^2$, C. R. Acad. Sci. Paris, t. 319, S\'{e}rie I, p. 665-668, 1994.

\bibitem{}  W. Sierpinski, ``Pythagorean Triangles'' , The Scripta Mathematica Studies, No. 9, Yeshiva University, New York, 1962, reprinted DOVER 2003.

\bibitem{}  M. Teigen and D. Hadwin; On Generating Pythagorean Triples, Amer. Math. Monthly 78,4 (1971), 378-379. 

\end{thebibliography}
\end{document}